
\input amstex \input amssym.def \input amssym.tex \input amsppt.sty \nologo


\magnification=\magstephalf

\catcode`\@=12
\voffset=2.5truepc
\hoffset=2.5truepc
\hsize=5.15truein
\vsize=8.50truein 
\def\halfhsize{2.5truein}

\tolerance=6000
\parskip=0pt
\parindent=18pt

\baselineskip=13.2pt

 \abovedisplayskip=12pt plus3pt minus2pt
 \belowdisplayskip=12pt plus3pt minus2pt
 \abovedisplayshortskip=12pt plus3pt minus2pt
 \belowdisplayshortskip=12pt plus3pt minus2pt

 \def\lskipamount{12pt}
 \def\lskip{\vskip\lskipamount plus3pt minus2pt}
 \def\lbreak{\par \ifdim\lastskip<\lskipamount
  \removelastskip \penalty-200 \lskip \fi}

 \def\lnobreak{\par \ifdim\lastskip<\lskipamount
  \removelastskip \penalty200 \lskip \fi}

\font\sectionfont=cmbx10 at 10pt
\font\autit=cmti12 
\font\titrm=cmbx12 at 14pt

\def\ded#1{\rightline{{\it #1}}\medskip}

\def\sec#1{\goodbreak\vskip 1.5truepc\centerline{\hbox {{\sectionfont #1}}}
\nobreak\vskip1truepc\noindent\nobreak}

\def\refs{\vskip 1.5pc{\centerline {\bf References}}\vskip 4pt \noindent}


\newif\iftitlepage\titlepagetrue
\newtoks\runningtitle \runningtitle={}
\newtoks\runningauthor \runningauthor={}

\def\rhead{\iftitlepage{\hfil\global\titlepagefalse}\else
{\rm\hfil{\it \the\runningtitle}\hfill\folio}\fi}
\def\lhead{\rm\folio\hfill{\global\titlepagefalse \it\the\runningauthor}}

\headline={\ifodd\pageno\rhead\else\lhead\fi}

\footline={\hfil}


\def\author{\bgroup\autit\bgroup}
\def\endauthor{\egroup\vskip 2truepc\egroup}

\def\title{\bgroup\bgroup\titrm\baselineskip20pt}
\def\endtitle{\egroup\vskip2truepc\egroup}

\def\Refs{\bgroup\refs\refstyle{A}\let\no=\key} \def\endRefs{\egroup}

\def\dedicatory{\bgroup\ded\bgroup}
\def\enddedicatory{\egroup\egroup}

\long\def\abstract#1\endabstract{
{\bf Abstract.} #1\par\vskip1.5pc}

\def\head#1\endhead{\sec{#1}}

\def\proclaim#1{\medskip{\noindent \bf  #1.}\bgroup\it}
\def\endproclaim{\egroup\medskip}

\def\demo#1{\medskip{\noindent\bf #1.}\rm\bgroup}
\def\enddemo{\egroup\medskip}

\def\remark#1{\medskip{\noindent\bf #1.}\rm\bgroup}
\def\endremark{\egroup\medskip}


\newbox\addressbox
\def\address{\setbox\addressbox=\vtop\bgroup
\hsize=\halfhsize\vglue3truepc\parindent=0pt
\obeylines}
\def\endaddress{\egroup} 

\newbox\addressboxtwo
\def\addresstwo{\setbox\addressboxtwo=\vtop\bgroup
\hsize=\halfhsize\vglue3truepc\parindent=0pt
\obeylines}
\def\endaddresstwo{\egroup}



\noindent SIGSAM Bulletin of the ACM
\vskip.01in
\noindent Vol. 30\#2, June 1996, Issue 116.

\bigskip

\title 
\centerline{Searching Symbolically for}
\centerline{Ap\'ery-like Formulae for Values of} 
\centerline{The Riemann Zeta Function}
\endtitle

\runningtitle{SYMBOLIC SEARCHING}
\runningauthor{BORWEIN \& BRADLEY}

\author 
\centerline{Jonathan Borwein and David Bradley\footnote"\dag"{Research supported by NSERC,
the Natural Sciences and Engineering Research Council of Canada.}}
\endauthor


\address 
Jonathan Borwein \& 
David Bradley,
Centre for Experimental
and Constructive Mathematics,
Simon Fraser University,
Burnaby B.C.,
Canada V5A 1S6.
jborwein\@cecm.sfu.ca
dbradley\@cecm.sfu.ca
\endaddress


\redefine\ge{\geqslant}
\redefine\i{\infty}

\define\({\left(}
\define\){\right)}
\define\[{\left[}
\define\]{\right]}

\define\C{\text{\bf C}}

\define\z{\zeta}
\redefine\l{\lambda}

\define\Li#1{\operatorname{Li_#1}}

\interdisplaylinepenalty=500

\abstract
We discuss some aspects of the search for identities
using computer algebra and symbolic methods.  To keep
the discussion as concrete as possible, we shall 
focus on so-called Ap\'ery-like formulae for special
values of the Riemann Zeta function.  Many of these results
are apparently new, and much more work needs to be done
before they can be formally proved and properly classified. 
A first step in this direction can be found in \cite{1}.
\endabstract

\head 1. Introduction \endhead
The Riemann Zeta function is 
$$
   \z(s) = \sum_{k=1}^\i \dfrac{1}{k^s},\qquad \Re(s)>1.\tag{1.1}
$$
In view of the ``Ap\'ery-like'' formulae
$$
   \z(2) = 3\sum_{k=1}^\i \dfrac{1}{k^2 {2k\choose k}},\quad
   \z(3) = \dfrac{5}{2}\sum_{k=1}^\i \dfrac{(-1)^{k+1}}{k^3 {2k\choose k}},\quad
   \z(4) = \dfrac{36}{17}\sum_{k=1}^\i \dfrac{1}{k^4 {2k\choose k}},\tag{1.2}
$$
one is tempted to speculate that there is an analogous formula for $\z(5)$, $\z(6)$,
$\z(7)$ and so on.  The key word here is {\it analogous}.  
For example, extensive computation has ruled out the possibility of formulae of the form
$$
   \z(5) = \dfrac{a}{b}\sum_{k=1}^\i \dfrac{(-1)^{k+1}}{k^5 {2k\choose k}},\quad
   \z(6) = \dfrac{c}{d}\sum_{k=1}^\i \dfrac{1}{k^6 {2k\choose k}},
$$
where $a$, $b$, $c$, $d$ are moderately sized integers.
Such negative results
are useful, as they tell us it would be a waste of time to search for interesting formulae
of a given form.  Thus, it would seem
there are no corresponding Ap\'ery-like formulae for higher zeta values.  End of story.
Consider however, the following result of 
Koecher \cite{2, 3}:
$$
   \z(5) = 2\sum_{k=1}^\i \dfrac{(-1)^{k+1}}{k^5{2k\choose k}} - \dfrac{5}{2}
            \sum_{k=1}^\i \dfrac{(-1)^{k+1}} {k^3{2k\choose k}} 
              \sum_{j=1}^{k-1}\dfrac{1}{j^2}.\tag{1.3}
$$
Koecher's formula points up a potential problem with symbolic searching.  Namely, negative
results need to be interpreted carefully, lest they be given more weight
than they deserve and unnecessarily discourage further investigation.  Also, it becomes clear
that symbolic searching is very much limited by the need to know fairly precisely the form
of what one is searching for in advance.

Koecher's formula (1.3) suggests that one might profit by searching for a formula of the form
$$
  \z(7) = r_1 \sum_{k=1}^\i \dfrac{(-1)^{k+1}}{k^7 {2k\choose k}} 
        + r_2 \sum_{k=1}^\i \dfrac{(-1)^{k+1}}{k^5 {2k\choose k}}\sum_{j=1}^{k-1}\dfrac{1}{j^2}
        + r_3 \sum_{k=1}^\i \dfrac{(-1)^{k+1}}{k^3 {2k\choose k}}\sum_{j=1}^{k-1}\dfrac{1}{j^4},
$$
where $r_1$, $r_2$, $r_3$ are rational numbers.  The following (conjectured)\footnote{See note below.} 
formula for $\z(7)$
was found \cite{1} using high precision arithmetic and Maple's integer relations algorithms:
$$
   \z(7) = \dfrac{5}{2} \sum_{k=1}^\i \dfrac{(-1)^{k+1}}{k^7 {2k\choose k}}
          +\dfrac{25}{2} \sum_{k=1}^\i \dfrac{(-1)^{k+1}}{k^3 {2k\choose k}}
                   \sum_{j=1}^{k-1} \dfrac{1}{j^4}.\tag{1.4}
$$
More generally, we have the (conjectured)\footnote{These conjectures have subsequently been 
proved.  See Granville and Almqvist's preprint 
http://www.math.uga.edu/~andrew/Postscript/BorBrad.ps.}
generating function formula \cite{1}
$$
   \sum_{k=1}^\i \dfrac{1}{k^3\(1-z^4/k^4\)}
 = \dfrac{5}{2}\sum_{k=1}^\i \dfrac{(-1)^{k+1}}{k^3 {2k\choose k}} \dfrac{1}{1-z^4/k^4}
                  \prod_{j=1}^{k-1}\dfrac{j^4+4z^4}{j^4-z^4},\quad z\in\C.\tag{1.5}
$$
Note that the constant coefficient in (1.5) gives the formula for $\z(3)$ in (1.2).  
The coefficient of $z^4$ in (1.5) gives (1.4).  We arrived at (1.5) by extensive use
of Maple's lattice algorithms, combined with a good deal of insightful guessing.  Interestingly,
Maple's convert(series, ratpoly) feature played a significant role.  The reader is referred to
\cite{1} for details.

Comparing our generating function
formula (1.5) with Koecher's \cite{3} 
$$
   \sum_{k=1}^\i \dfrac{1}{k^3\(1-z^2/k^2\)} 
 = \sum_{k=1}^\i \dfrac{(-1)^{k+1}}{k^3 {2k\choose k}}\(\dfrac{1}{2}+\dfrac{2}{1-z^2/k^2}\)
                   \prod_{j=1}^{k-1}(1-z^2/j^2)\tag{1.6}
$$
raises some interesting issues related to formula redundancy, and which remain unresolved.
We address certain of these issues in the next section.


\head 2. Redundancy Relations\endhead
To mitigate the problem of symbol clutter in what follows requires some notation.  
We denote the power sum symmetric functions by
$$
   P_r(k) := \left\{\aligned \sum_{j=1}^{k-1} j^{-r},\quad r\ge 1,\\
                                   1,\quad r=0.
                          \endaligned\right.
$$
Next, we define functions $\l$, $\mu$ by
$$
\align
   \l(m,\prod_{j=1}^n P_{r_j}) 
&:= \sum_{k=1}^\i \dfrac{(-1)^{k+1}}{k^m {2k\choose k}}\prod_{j=1}^n P_{r_j}(k),\\
   \mu(m,\prod_{j=1}^n P_{r_j}) 
&:= \sum_{k=1}^\i \dfrac{1}{k^m {2k\choose k}}\prod_{j=1}^n P_{r_j}(k).
\endalign
$$
In the new notation, (1.2) becomes
$$
   \z(2)=3\mu(2,P_0),\qquad \z(3) = \dfrac{5}{2}\l(3,P_0),\qquad \z(4) = \dfrac{36}{17}\mu(4,P_0),\tag{2.1}
$$
while (1.3) and (1.4) become
$$
   \z(5) = 2\l(5,P_0) -\dfrac{5}{2}\l(3,P_2),\qquad 
   \z(7) = \dfrac{5}{2}\l(7,P_0)+\dfrac{25}{2}\l(3,P_4),\tag{2.2}
$$
respectively.
To illustrate the issue of formula redundancy, consider Koecher's formula for $\z(7)$ \cite{3}
which becomes, in our notation,
$$
  \z(7) 
= 2\l(7,P_0) - 2\l(5,P_2) + \dfrac{5}{4}\l(3,P_2^2) - \dfrac{5}{4}\l(3,P_4).\tag{2.3}
$$
In view of the second formula in (2.2),
the middle two terms of (2.3) must be redundant.  Indeed, lattice-based reduction shows that
$$
  -2\l(5,P_2) + \dfrac{5}{4}\l(3,P_2^2) 
= \dfrac{55}{4}\l(3,P_4)+\dfrac{1}{2}\l(7,P_0).\tag{2.4}
$$
Although we currently have no real understanding why interrelations between $\l$ sums such as (2.4)
hold, we decided to limit our symbolic search for Zeta function identities in which no such interrelations
exist.\footnote{Of course, we cannot prove that (1.4) contains no redundancy, since, for example, we 
cannot even prove that $\z(7)$ is irrational.}  This was carried out by starting with a ``full set'' of $\l$
sums and checking that a relation holds with the relevant Zeta value.  Now recurse, using the following scheme.
From any found relation, toss out the Zeta value.
If no relation is found amongst the remaining sums, output the relation that held when the Zeta value
was included, and report it as non-redundant.
Otherwise, systematically 
discard the various $\l$ sums from the list,
until a non-redundant relation remains.  Carrying out the aforementioned procedure yields the following
formulae which evidently exhaust the list of non-redundant formulae for each given Zeta 
value:
$$
\align
 & 17\z(4)-36\mu(4,P_0) = 5\z(4)-108\mu(2,P_2) = 0,\\
 &\\
 &7\z(6)+1944\mu(2,P_4)-1944\mu(2,P_2^2)\\
&= 215\z(6)-2592\mu(4,P_2)-3888\mu(2,P_4)\\
&=229\z(6)-2592\mu(4,P_2)-3888\mu(2,P_2^2)\\
&=1481\z(6)-2592\mu(6,P_0)-3888\mu(2,P_2^2)\\
&=313\z(6)-648\mu(6,P_0)+648\mu(4,P_2)\\
&=163\z(6)-288\mu(6,P_0)-432\mu(2,P_4)\\
&=0,\\
&\\
 & 2\z(7)-5\l(7,P_0)-25\l(3,P_4) \\
&= 4\z(7)-25\l(3,P_2^2)+40\l(5,P_2)+225\l(3,P_4)\\
&= 22\z(7) - 25\l(3,P_2^2)+40\l(5,P_2)-45\l(7,P_0)\\
&= 0,\\
&\\
 & 72\z(9)+135\l(7,P_2)-147\l(9,P_0)-60\l(5,P_2^2)-85\l(3,P_6)+25\l(3,P_2^3)\\
&= 36\z(9) -540\l(5,P_4)-96\l(9,P_0)+60\l(5,P_2^2)-1130\l(3,P_6)\\
&\quad+675\l(3,P_4P_2)-25\l(3,P_2^3)\\
&= 4\z(9)+196\l(5,P_4)+32\l(7,P_2)-36\l(5,P_2^2)+390\l(3,P_6)\\
&\quad-245\l(3,P_4P_2)+15\l(3,P_2^3)\\
&= 4\z(9)-20\l(5,P_4)+5\l(7,P_2)-9\l(9,P_0)-45\l(3,P_6)+25\l(3,P_4P_2)\\ 
&= 116\z(9)+68\l(5,P_4)+226\l(7,P_2)-234\l(9,P_0)-108\l(5,P_2^2)\\
&\quad-85\l(3,P_4P_2)+45\l(3,P_2^3)\\
&=0.
\endalign
$$

No additional formulae other than the formulae given in \S 1 were found for
$\z(2)$, $\z(3)$ and $\z(5)$.  We discuss additional uniqueness issues in the
next section.  

\head 3. Uniqueness and $\z(4n+3)$ \endhead
If one extends the list given in the previous section, it becomes
apparent that $\z(4n+3)$ evidently has a unique representation in terms of
$\l$ sums of the form $\l(m,P_r)$ in which $r$ is always a multiple of
four.  We exploited this observation in \cite{1} to arrive at our generating
function formula (1.5).  Unfortunately, there seems to be no sensible selection
to make amongst the formulae for $\z(4n+1)$ which gives an analogous generating
function identity.  Our 
Maple code for producing all
possible non-redundant formulae for $\z(13)$
ran for over two months before it was killed.
The resulting incomplete file is over three
thousand lines long and contains hundreds and hundreds of independent formulae.
If a generating function identity (other than a bisection of Koecher's)
for $\z(4n+1)$ is found, it is unlikely that
it will be discovered by hunting for the appropriate representatives from the
identities for $\z(9)$, $\z(13)$, etc. and looking for a pattern.  

Recall Ramanujan's formulae \cite{4}
$$
\align
   2\z(4n+3) 
&= (2\pi)^{4n+3}\sum_{k=0}^{2n+2} (-1)^{k+1} \dfrac{B_{2k}}{(2k)!} 
                                             \dfrac{B_{4n+4-2k}}{(4n+4-2k)!}
   -4\sum_{k=1}^\i \dfrac{k^{-4n-3}}{e^{2\pi k}-1}\\
\endalign
$$
and
$$
\align
   2\z(4n+1)
&= (2\pi)^{4n+1} \dfrac{1}{2n}\sum_{k=0}^{2n+1} (-1)^{k+1} (2k-1) \dfrac{B_{2k}}{(2k)!}
                         \dfrac{B_{4n+2-2k}}{(4n+2-2k)!}\\
&\quad  -4\sum_{k=1}^\i \dfrac{k^{-4n-1}}{e^{2\pi k}-1}
   -\dfrac{\pi}{n}\sum_{k=1}^\i \dfrac{k^{-4n}}{\sinh^2(\pi k)}.
\endalign
$$
Here, the additional complexity in the $4n+1$ case arises from taking the derivative of the
appropriate modular transformation formula.
Perhaps there is an analogous phenomenon operating in the case of these Ap\'ery-like
identities as well.

\copy\addressbox
\bye